%% file: syncct.tex
\documentclass{article}  
\input{macros}

\begin{document}
\title{Synchronizing continuous-time neutrally stable linear systems via 
partial-state coupling} 
\author{S. Emre Tuna\\
{\small {\tt tuna@eee.metu.edu.tr}} }
\maketitle

\begin{abstract}                          
Synchronization of coupled continuous-time linear systems is studied
in a general setting. For identical neutrally-stable linear systems
that are detectable from their outputs, it is shown that a linear
output feedback law exists under which the coupled systems globally
asymptotically synchronize under all fixed (directed) connected network
topologies. An algorithm is provided to compute one such feedback law
based on individual system parameters. The dual case, where individual
systems are neutrally stable and stabilizable from their inputs, is
also considered and parallel results are established.
\end{abstract}

\section{Introduction}

In \cite{aut7479} we have shown, for identical discrete-time linear
systems that are detectable (stabilizable) from their outputs (inputs)
and neutrally stable, that a linear feedback law exists under which
the coupled systems globally asymptotically synchronize for all fixed
(directed) connected network topologies. There we have also provided
an algorithm to compute such feedback law based on individual system
parameters. In this companion paper we provide counterpart results for
continuous-time linear systems.

\subsection{Background}

``The main issue in studying the synchronization of coupled dynamical
systems is the stability of synchronization. As in all cases where
stability is the issue, the question whose answer is sought is {\em
Under what conditions} will the individual systems synchronize? In a
simplified yet widely-studied scenario, where the individual system
dynamics are identical and the coupling between them is linear,
studies focus on two ingredients: the dynamics of an individual system
and the network topology. Starting with the {\em agreement algorithm}
in \cite{tsitsiklis86} a number of contributions
\cite{jadbabaie03,moreau05,ren05,angeli06,olfati04} have
gathered around the case where the weakest possible assumptions are
made on the network topology at the expense of restrictive individual
system dynamics. It was established in those works on {\em multi-agent
systems} \cite{olfati07} that when the individual system is taken to
be an integrator and the coupling is of full-state, synchronization
({\em consensus}) results for time-varying interconnections whose
unions\footnote{By {\em union of interconnections} we actually mean
the union of the graphs representing the interconnections.} over an
interval are assumed to be connected instead of that each
interconnection at every instant is connected.
         
``Another school of research investigates networks with more complicated
(nonlinear) individual system dynamics \cite{strogatz01,wang02}. When
that is the case, the restrictions on the network topology have to be
made stricter in order to ensure stability of
synchronization. Generally speaking, more than mere connectedness of
the network has been needed: coupling strength is required to be
larger than some threshold and sometimes a symmetry or balancedness
assumption is made on the connection graph. Different (though related)
approaches have provided different insights over the years. The
primary of such approaches is based on the calculations of the
eigenvalues of the connection matrix and a parameter (e.g. the maximal
Lyapunov exponent) depending on the individual system dynamics
\cite{wu95,pecora98,chen03}. In endeavor to better understand
synchronization stability, tools from systems theory such as Lyapunov
functions \cite{belykh06,hui07}, passivity
\cite{pogromsky01,arcak07,stan07,wu01}, contraction theory
\cite{slotine04}, and incremental input-to-state stability
($\delta$ISS) theory \cite{cai06} have also proved
useful.''\footnote{Borrowed from \cite{aut7479}.}

\subsection{Contribution}

In this paper we study two dual problems. In the first case we
consider the following individual system 
\begin{eqnarray}\label{eqn:C}
\dot{x}_{i}=Ax_{i}\,,\quad y_{i}=Cx_{i}\,,
\end{eqnarray}
where $A$ is assumed to be {\em neutrally stable} and pair $(C,\,A)$
{\em detectable}, and design a linear output feedback
gain $L$ that synchronizes any fixed {\em connected} network of any number
of coupled replicas of \eqref{eqn:C}. Such $L$ guarantees the
synchronization of $p$ individual systems when coupled as
\begin{eqnarray*}
\dot{x}_{i}=Ax_{i}+L\sum_{j=1}^{p}\gamma_{ij}(y_{j}-y_{i})\,.
\end{eqnarray*}
As the dual problem we consider 
\begin{eqnarray}\label{eqn:B}
\dot{x}_{i}=Ax_{i}+Bu_{i}\,,
\end{eqnarray}
where $A$ is assumed to be {\em neutrally stable} and pair $(A,\,B)$  
{\em stabilizable}, and design a linear feedback gain $K$ that 
synchronizes any fixed {\em connected} network of any number
of coupled replicas of \eqref{eqn:B}. Such $K$ guarantees the
synchronization of $p$ individual systems when coupled as
\begin{eqnarray*}
\dot{x}_{i}=Ax_{i}+BK\sum_{j=1}^{p}\gamma_{ij}(x_{j}-x_{i})\,.
\end{eqnarray*}
To the best of our knowledge, feedback design (in such a general
setting) in order to guarantee synchronization under arbitrary (fixed)
interconnections is a novelty of our work. It is worth noting that
our main theorems make a compromise result between the two previously
mentioned cases (i) where synchronization is established for very
primitive individual system dynamics, such as that of an integrator,
but under the weakest conditions on the network topology and (ii)
where the network topology has to satisfy stronger conditions, such as
that the coupling strength should be above a threshold, for want of
achieving synchronization for nonlinear individual system dynamics.

\subsection{Organization}
The remainder of the paper is organized as follows. In the next
section we provide notation and some preliminaries.  Then we formally
state our problems in Sections~\ref{sec:problemstatement} and
\ref{sec:dualproblem}.  Section~\ref{sec:pre} is where we establish
our key result which we will later use to solve the problems we aim
at. In Section~\ref{sec:main} we provide an algorithm to design output
feedback gain that we seek for synchronization and prove that it
works. Then, in Section~\ref{sec:maindual}, we design a state feedback
gain that solves the dual problem.

\section{Notation and definitions}

Let $\Natural$ denote the set of nonnegative integers and $\Real_{\geq
0}$ set of nonnegative real numbers. Let $|\cdot|$ denote 2-norm.
Identity matrix in $\Real^{n\times n}$ is denoted by $I_{n}$. A matrix
$A\in\Real^{n\times n}$ is {\em Hurwitz} if all of its eigenvalues
have strictly negative real parts. A matrix $S\in\Real^{n\times n}$ is
{\em skew-symmetric} if $S+S^{T}=0$. Given $C\in\Real^{m\times n}$ and
$A\in\Real^{n\times n}$, pair $(C,\,A)$ is {\em observable} if
$[C^{T}\, A^{T}C^{T}\, A^{2T}C^{T}\,\ldots\, A^{(n-1)T}C^{T}]$ is full
row rank. Pair $(C,\,A)$ is {\em detectable} (in the continuous-time
sense) if that $Ce^{At}x=0$ for some $x\in\Real^{n}$ and for all
$t\geq 0$ implies $\lim_{t\to\infty}e^{At}x=0$. Given
$B\in\Real^{n\times m}$ and $A\in\Real^{n\times n}$, pair $(A,\,B)$ is
{\em controllable (stabilizable)} if $(B^{T},\,A^{T})$ is observable
(detectable). Matrix $A\in\Real^{n\times n}$ is {\em neutrally stable}
(in the continuous-time sense) if it has no eigenvalue with positive
real part and the Jordan block corresponding to any eigenvalue on the
imaginary axis is of size one.\footnote{Note that $A$ is neutrally
stable iff there exists a symmetric positive definite matrix $P$ such
that $A^{T}P+P{A}\leq 0$,
\cite{antsaklis97}.}  Let $\one\in\Real^{p}$ denote the vector with
all entries equal to one.

{\em Kronecker product} of $A\in\Real^{m\times n}$ and $B\in\Real^{p\times q}$ is
\begin{eqnarray*}
A\otimes B:=
\left[
\begin{array}{ccc}
a_{11}B & \cdots & a_{1n}B\\
\vdots  & \ddots & \vdots\\
a_{m1}B & \cdots & a_{mn}B
\end{array}
\right]
\end{eqnarray*}
Kronecker product comes with the properties $(A\otimes
B)(C\otimes D)=(AC)\otimes(BD)$ (provided that products $AC$ and $BD$
are allowed) $A\otimes B+A\otimes C=A\otimes(B+C)$ (for $B$ and
$C$ that are of same size) and $(A\otimes B)^{T}=A^{T}\otimes B^{T}$. 

A ({\em directed}) {\em graph} is a pair $(\N,\,\A)$ where $\N$ is a
nonempty finite set (of {\em nodes}) and $\A$ is a finite collection
of pairs ({\em arcs}) $(n_{i},\,n_{j})$ with $n_{i},\,n_{j}\in\N$. A
{\em path} from $n_{1}$ to $n_{\ell}$ is a sequence of nodes
$\{n_{1},\,n_{2},\,\ldots,\,n_{\ell}\}$ such that $(n_{i},\,n_{i+1})$
is an arc for $i\in\{1,\,2,\,\ldots,\,\ell-1\}$. A graph is {\em
connected} if it has a node to which there exists a path from every
other node.\footnote{Note that this definition of connectedness for
directed graphs is weaker than strong connectivity and stronger than
weak connectivity.}

The graph of a matrix $\Gamma:=[\gamma_{ij}]\in\Real^{p\times p}$ is
the pair $(\N,\,\A)$ where $\N =\{n_{1},\,n_{2},\,\ldots,\,n_{p}\}$
and $(n_{i},\,n_{j})\in\A$ iff $\gamma_{ij}>0$. Matrix $\Gamma$ is
said to be {\em connected} (in the continuous-time sense) if it
satisfies:
\begin{enumerate}
\item[(i)] $\gamma_{ij}\geq 0$ for $i\neq j$;
\item[(ii)] each row sum equals 0;
\item[(iii)] its graph is connected.
\end{enumerate}

For connected $\Gamma$, it follows from definition that $\lambda=0$ is
an eigenvalue with eigenvector $\one$ (i.e. $\Gamma\one=0$.) Moreover,
all the other eigenvalues have real parts strictly negative. Let
$r^{T}$ be the left eigenvector of eigenvalue $\lambda=0$
(i.e. $r^{T}\Gamma=0$) with $r^{T}\one=1$. Then
$\lim_{t\to\infty}e^{\Gamma{t}}={\one}r^{T}$.

Given maps $\xi_{i}:\Real_{\geq 0}\to\Real^{n}$ for
$i=1,\,2,\,\ldots,\,p$ and a map $\bar\xi:\Real_{\geq
0}\to\Real^{n}$, the elements of the set
$\{\xi_{i}(\cdot):i=1,\,2,\,\ldots,\,p\}$ are said to {\em synchronize
to} $\bar{\xi}(\cdot)$ if $|\xi_{i}(t)-\bar\xi(t)|\to 0$ as
$t\to\infty$ for all $i$.

\section{Problem I}\label{sec:problemstatement}  
We now formalize our first problem.

\subsection{Systems under study}
We consider $p$ identical linear systems
\begin{eqnarray}\label{eqn:system}
{\dot{x}}_{i}=Ax_{i}+u_{i}\, ,\quad y_{i}=Cx_{i}\, ,\quad i=1,\,2,\,\ldots,\,p
\end{eqnarray}
where $x_{i}\in\Real^{n}$ is the {\em state}, $u_{i}\in\Real^{n}$ is the {\em input}, and
$y_{i}\in\Real^{m}$ is the {\em output} of the $i$th system. Matrices
$A$ and $C$ are of proper dimensions. The solution of $i$th system at
time $t\geq 0$ is denoted by $x_{i}(t)$. In this paper we consider
the case where at each time instant only the following information
\begin{eqnarray}\label{eqn:z}
z_{i}&=&\sum_{j=1}^{p}\gamma_{ij}(y_{j}-y_{i})
\end{eqnarray}
is available to $i$th system to determine an input value where
$\gamma_{ij}$ are the entries of the matrix $\Gamma\in\Real^{p\times p}$ describing
the network topology. Nondiagonal entries of $\Gamma$ are nonnegative
and each row sums up to zero. That is, the coupling between systems is
diffusive.

\subsection{Assumptions made}
We make the following assumptions on systems~\eqref{eqn:system} which
will henceforth hold.\\

\noindent
{\bf (A1)} $A$ is neutrally stable.\\
{\bf (A2)} $(C,\,A)$ is detectable.

\subsection{Objectives}
Our first objective is to {\em show that there exists a
linear feedback law $L\in\Real^{n\times m}$ such that, for all $p$ and
connected $\Gamma\in\Real^{p\times p}$, solutions of
systems~\eqref{eqn:system} with $u_{i}=Lz_{i}$, where $z_{i}$ is as in
\eqref{eqn:z}, globally (i.e. for all initial conditions)
synchronize to a bounded trajectory.} Our second objective is to {\em devise an algorithm to
compute one such $L$.}

\section{Problem II}\label{sec:dualproblem}  
In this section we state the second problem, which, as noted earlier,
is the dual of the first.
\subsection{Systems under study}
Consider $p$ identical linear systems
\begin{eqnarray}\label{eqn:dualsystem}
{\dot{x}}_{i}=Ax_{i}+Bu_{i}\, ,\quad i=1,\,2,\,\ldots,\,p
\end{eqnarray}
where $x_{i}\in\Real^{n}$ and $u_{i}\in\Real^{m}$. Matrices
$A$ and $B$ are of proper dimensions. We consider
the case where at each time instant the following information
\begin{eqnarray}\label{eqn:zdual}
z_{i}&=&\sum_{j=1}^{p}\gamma_{ij}(x_{j}-x_{i})
\end{eqnarray}
is available to $i$th system to determine an input value.

\subsection{Assumptions made}
We make the following assumptions on systems~\eqref{eqn:dualsystem}
which will henceforth hold.\\

\noindent
{\bf (B1)} $A$ is neutrally stable.\\
{\bf (B2)} $(A,\,B)$ is stabilizable.

\subsection{Objectives}
Our first objective regarding the dual problem is to {\em show that
there exists a linear feedback law $K\in\Real^{m\times n}$ such that,
for all $p$ and connected $\Gamma\in\Real^{p\times p}$, solutions of
systems~\eqref{eqn:dualsystem} with $u_{i}=Kz_{i}$, where $z_{i}$ is
as in
\eqref{eqn:zdual}, globally (i.e. for all initial conditions)
synchronize to a bounded trajectory.} Our second objective is to {\em
devise an algorithm to compute one such $K$.}

\section{A special case}\label{sec:pre}
Before we attempt to solve Problems I and II, we first establish a
preliminary result to be resorted later. Consider the following
coupled systems
\begin{eqnarray}\label{eqn:systemct}
\dot\xi_{i}=S\xi_{i}+H^{T}H\sum_{j=1}^{p}\gamma_{ij}(\xi_{j}-\xi_{i})\, ,\quad i=1,\,2,\,\ldots,\,p
\end{eqnarray}
where $\xi_{i}\in\Real^{n}$ is the state of the $i$th system,
$S\in\Real^{n\times n}$, and $H\in\Real^{m\times n}$. We make the
following assumptions on systems~\eqref{eqn:systemct} which will
henceforth hold.\\

\noindent
{\bf (C1)}\ $S$ is skew-symmetric.\\
{\bf (C2)}\ $(H,\,S)$ is observable.\\
{\bf (C3)}\ $\Gamma:=[\gamma_{ij}]$ is connected.\\

\noindent
Below we provide our first result.

\begin{theorem}\label{thm:lyap}
Consider systems~\eqref{eqn:systemct}. Let $r\in\Real^{p}$ be such
that $r^{T}\Gamma=0$ and $r^{T}\one=1$. Then solutions
$\xi_{i}(\cdot)$, for $i=1,\,2,\,\ldots,\,p$, synchronize to
\begin{eqnarray*}
\bar{\xi}(t):=({r^{T}\otimes e^{St}})
\left[\!\!
\begin{array}{c}
\xi_{1}(0)\\
\vdots\\
\xi_{p}(0)
\end{array}
\!\!
\right]
\end{eqnarray*}
\end{theorem}

\begin{proof}
Consider matrix $\Gamma-\one{r}^{T}$. Observe that
$(\Gamma-\one{r}^{T})^{k}=\Gamma^{k}+(-1)^{k}\one{r^{T}}$ for $k\in\Natural$. For
$t\in\Real$ therefore we can write
\begin{eqnarray*}
e^{(\Gamma-1r^{T})t}
&=&I_{p}+t(\Gamma-1r^{T})+\frac{t^{2}}{2}(\Gamma-1r^{T})^{2}+\ldots\\
&=&\left(I_{p}+t\Gamma+\frac{t^{2}}{2}\Gamma^{2}+\ldots\right)
-\left(t\one{r^{T}}-\frac{t^{2}}{2}\one{r}^{T}+\ldots\right)\\
&=&e^{\Gamma{t}}-(1-e^{-t})\one{r^{T}}\,.
\end{eqnarray*} 
Consequently $\lim_{t\to\infty}e^{(\Gamma-\one{r^{T}})t}=0$. We deduce
therefore that $\Gamma-\one{r^{T}}$ is Hurwitz. Since
$\Gamma-\one{r^{T}}$ is Hurwitz, there exist symmetric positive
definite matrices $P,\,Q\in\Real^{p\times p}$ such that
\begin{eqnarray}\label{eqn:prepost}
-Q=(\Gamma-\one{r^{T}})^{T}P+P(\Gamma-\one{r^{T}})\,.
\end{eqnarray}
Define positive semidefinite matrices
$\widehat{P}:=(I_{p}-\one{r^{T}})^{T}P(I_{p}-\one{r^{T}})$ and
$\widehat{Q}:=(I_{p}-\one{r^{T}})^{T}Q(I_{p}-\one{r^{T}})$. Now pre-
and post-multiply equation \eqref{eqn:prepost} by
$(I_{p}-\one{r^{T}})^{T}$ and $(I_{p}-\one{r^{T}})$, respectively. We
obtain
\begin{eqnarray*}
-\widehat{Q}
&=&(I_{p}-\one{r^{T}})^{T}(\Gamma-\one{r^{T}})^{T}P(I_{p}-\one{r^{T}})\nonumber\\
&&\qquad+(I_{p}-\one{r^{T}})^{T}P(\Gamma-\one{r^{T}})(I_{p}-\one{r^{T}})\nonumber\\
&=&\Gamma^{T}P(I_{p}-\one{r^{T}})+(I_{p}-\one{r^{T}})^{T}P\Gamma\nonumber\\
&=&\Gamma^{T}(I_{p}-\one{r^{T}})^{T}P(I_{p}-\one{r^{T}})
+(I_{p}-\one{r^{T}})^{T}P(I_{p}-\one{r^{T}})\Gamma\nonumber\\
&=&\Gamma^{T}\widehat{P}+\widehat{P}\Gamma\,.
\end{eqnarray*}
We now stack the individual system states to obtain
$\ex:=[\xi_{1}^{T}\ \xi_{2}^{T}\ \cdots\ \xi_{p}^T]^{T}$. We can then cast
\eqref{eqn:systemct} into
\begin{eqnarray}\label{eqn:wrt}
\dot{\ex}=(I_{p}\otimes{S}+\Gamma\otimes{H^{T}H})\ex\,.
\end{eqnarray} 
Define $V:\Real^{pn}\to\Real_{\geq 0}$ as $V(\ex):=\ex^{T}(\widehat{P}\otimes I_{n})\ex$. 
Differentiating $V(\ex(t))$ with respect to time we obtain
\begin{eqnarray}\label{eqn:vdot}
\dot{V}(\ex)
&=&\ex^{T}(I_{p}\otimes{S^{T}}+\Gamma^{T}\otimes{H^{T}H})(\widehat{P}\otimes I_{n})\ex\nonumber\\
&&\qquad+\ex^{T}(\widehat{P}\otimes I_{n})(I_{p}\otimes{S}+\Gamma\otimes{H^{T}H})\ex\nonumber\\
&=&\ex^{T}(\widehat{P}\otimes({S^{T}}+S)
+(\Gamma^{T}\widehat{P}+\widehat{P}\Gamma)\otimes{H^{T}H})\ex\nonumber\\
&=&-\ex^{T}(\widehat{Q}\otimes{H^{T}H})\ex\,.
\end{eqnarray}
Thence $\dot{V}(\ex)\leq 0$ for both $\widehat{Q}$ and $H^{T}H$ (and
consequently their Kronecker product) are positive semidefinite.

Given some $\zeta\in\Real^{pn}$, let $\X\subset\Real^{pn}$ be the
closure of the set of all points $\eta$ such that
$\eta=(\one{r^{T}}\otimes e^{St})\zeta$ for some $t\geq 0$. Set $\X$
is compact for it is closed by definition and bounded due to that
$\zeta$ is fixed and $S$ is a neutrally-stable matrix. Having defined
$\X$, we now define 
\begin{eqnarray*}
\Omega:=\{\eta\in\Real^{pn}:(\one{r^{T}}\otimes
I_{n})\eta\in\X\,,V(\eta)\leq V(\zeta)\}\,.  
\end{eqnarray*}
Let us show that $\Omega$ is forward invariant. Observe that
\begin{eqnarray*}
\frac{d}{dt}\left((\one{r^{T}}\otimes I_{n})\ex(t)\right)
&=&(\one{r^{T}}\otimes I_{n})(I_{p}\otimes{S}+\Gamma\otimes{H^{T}H})\ex(t)\\
&=&(\one{r^{T}}\otimes S + \one{r^{T}}\Gamma\otimes{H^{T}H})\ex(t)\\
&=&(\one{r^{T}}\otimes S)\ex(t)\\
&=&(I_{p}\otimes{S})(\one{r^{T}}\otimes I_{n})\ex(t)\,.
\end{eqnarray*} 
We therefore have 
\begin{eqnarray}\label{eqn:aksuleyman}
(\one{r^{T}}\otimes I_{n})\ex(t)=(\one{r^{T}}\otimes e^{St})\ex(0)
\end{eqnarray}
which in turn implies that if $(\one{r^{T}}\otimes I_{n})\ex(0)\in\X$
then $(\one{r^{T}}\otimes I_{n})\ex(t)\in\X$ for all $t\geq
0$. Likewise, if $V(\ex(0))\leq V(\zeta)$ then $V(\ex(t))\leq
V(\zeta)$ for all $t\geq 0$ thanks to \eqref{eqn:vdot}. As a result,
if $\ex(0)\in\Omega$ then $\ex(t)\in\Omega$ for all $t\geq 0$, that
is, $\Omega$ is forward invariant with respect to \eqref{eqn:wrt}.

Set $\Omega$ is closed by construction. To show that it is compact
therefore all we need to do is to establish its boundedness. Let 
\begin{eqnarray*}
a:=\sup_{V(\eta)\leq V(\zeta)}|\eta-(\one{r^{T}}\otimes I_{n})\eta|\,.
\end{eqnarray*}
If we go back to the definition of $V$ we immediately see that $a<\infty$. Now let
\begin{eqnarray*}
b:=\sup_{\omega\in\X}|\omega|\,.
\end{eqnarray*}
Since $\X$ is bounded, $b<\infty$ as well. Now, given any $\eta\in\Omega$ we have
$|\eta-(\one{r^{T}}\otimes I_{n})\eta|\leq a$. Hence we can write
\begin{eqnarray*}
|\eta|
&\leq& a+|(\one{r^{T}}\otimes I_{n})\eta|\\
&\leq& a+\sup_{\omega\in\X}|\omega|\\
&=& a+b\,.
\end{eqnarray*}
Therefore $\Omega$ is bounded. Having shown that $\Omega$ is forward
invariant and compact, we can now invoke LaSalle's invariance
principle \cite[Thm.~3.4]{khalil96} and claim that any solution
starting in $\Omega$ approaches to the largest invariant set
$\W\subset\{\eta\in\Omega:\dot{V}(\eta)=0\}$.

Let now $\eta(\cdot)$ be a solution of \eqref{eqn:wrt} such that $\eta(t)\in\W$ for all $t\geq 0$. 
Given some $\tau\geq 0$, since $\dot{V}(\eta(\tau))=0$, we can write 
\begin{eqnarray*}
0
&=&\eta(\tau)^{T}(\widehat{Q}\otimes H^{T}H)\eta(\tau)\\
&=&\eta(\tau)^{T}((I_{p}-\one{r^{T}})^{T}Q(I_{p}-\one{r^{T}})\otimes H^{T}H)\eta(\tau)
\end{eqnarray*} 
which implies, since $Q$ is positive definite, that either
$((I_{p}-\one{r^{T}})\otimes I_{n})\eta(\tau)=0$ or $(I_{p}\otimes
H)\eta(\tau)=0$. Suppose now that 
\begin{eqnarray}\label{eqn:star}
((I_{p}-\one{r^{T}})\otimes I_{n})\eta(\tau)\neq 0\,. 
\end{eqnarray}
Continuity of $\eta(\cdot)$ implies that there exists $\delta>0$ such
that $((I_{p}-\one{r^{T}})\otimes I_{n})\eta(t)\neq 0$ for
$t\in[\tau,\,\tau+\delta]$. Therefore we must have $(I_{p}\otimes
H)\eta(t)=0$ for $t\in[\tau,\,\tau+\delta]$. However, observability of
pair $(H,\,S)$ stipulates that $\eta(t)=0$ for
$t\in[\tau,\,\tau+\delta]$ which contradicts \eqref{eqn:star}. We
then deduce $((I_{p}-\one{r^{T}})\otimes I_{n})\eta(t)=0$ for all
$t\geq 0$.  Therefore
$\W\subset\{\omega\in\Omega:\omega=(\one{r^{T}}\otimes
I_{n})\omega\}=\X$.

Let us now be given any solution $\ex(\cdot)$ of \eqref{eqn:wrt}.
Since $\zeta$ that we used to construct $\Omega$ was arbitrary, without
loss of generality, we can take $\ex(0)=\zeta$. That $\ex(0)\in\Omega$
implies that $\ex(t)$ approaches $\X$ as $t\to\infty$. Therefore we
are allowed to write
\begin{eqnarray*}
0
&=&\lim_{t\to \infty}\left(\ex(t)-(\one{r^{T}}\otimes I_{n})\ex(t)\right)\\
&=&\lim_{t\to \infty}\left(\ex(t)-(\one{r^{T}}\otimes e^{St})\ex(0)\right)
\end{eqnarray*}
where we used \eqref{eqn:aksuleyman}.
\end{proof}
\\

\noindent
The following result (cf. \cite{hui07}) comes as a byproduct
of Theorem~\ref{thm:lyap}.

\begin{corollary}
Consider coupled harmonic oscillators (in $\Real^{2}$) described by
\begin{eqnarray*}
\dot{x}_{i}&=& y_{i}\\
\dot{y}_{i}&=& -x_{i}+\sum_{j=1}^{p}\gamma_{ij}(y_{j}-y_{i})\, ,\quad i=1,\,2,\,\ldots,\,p\,.
\end{eqnarray*}
Oscillators synchronize for all connected $\Gamma$.
\end{corollary}

\section{Solution to Problem I}\label{sec:main}

In this section we use Theorem~\ref{thm:lyap} in order to reach our objectives 
stated in Section~\ref{sec:problemstatement}. We first give the following 
fact.

\begin{fact}\label{fact:one}
Let $F\in\Real^{n\times n}$ be a neutrally-stable matrix with all its eigenvalues 
residing on the imaginary axis. Then 
\begin{eqnarray}\label{eqn:limit}
P:=\lim_{t\to\infty}\ t^{-1}\int_{0}^{t}e^{F^{T}\tau}e^{F\tau}d\tau
\end{eqnarray}
is well-defined and symmetric positive definite. It also satisfies $PF+F^{T}P=0$.
\end{fact}

\begin{proof}
Matrix $F$ is similar to a skew-symmetric matrix. Therefore $e^{Ft}$
is (almost) periodic \cite{vanvleck64}. Periodicity directly yields
that limit in \eqref{eqn:limit} exists, that is, $P$ is well-defined.
Similarity to a skew-symmetric matrix also brings that $\inf_{t\in\Real}|e^{Ft}|>0$
and $\sup_{t\in\Real}|e^{Ft}|<\infty$. Same goes for
$F^{T}$. Therefore there exist scalars $a,\,b>0$ such that $aI_{n}\leq
e^{F^{T}t}e^{Ft}\leq bI_{n}$ for all $t\in\Real$. We can then write
\begin{eqnarray*}
aI_{n}\leq t^{-1}\int_{0}^{t}e^{F^{T}\tau}e^{F\tau}d\tau\leq bI_{n}
\end{eqnarray*}
for all $t\geq 0$. Therefore $P$ is positive definite. Symmetricity of $P$ comes 
by construction. Finally, observe that
\begin{eqnarray*}
|PF+F^{T}P|
&=&\lim_{t\to\infty}\ t^{-1}\left|\int_{0}^{t}\left(e^{F^{T}\tau}e^{F\tau}F
+F^{T}e^{F^{T}\tau}e^{F\tau}\right)d\tau\right|\\
&=&\lim_{t\to\infty}\ t^{-1}\left|\int_{0}^{t}d\left(e^{F^{T}\tau}e^{F\tau}\right)\right|\\
&\leq&\lim_{t\to\infty}\ t^{-1}\left(\left|e^{F^{T}t}e^{Ft}\right|
+\left|e^{F^{T}0}e^{F0}\right|\right)\\
&\leq&\lim_{t\to\infty}\ t^{-1}(b+1)\\
&=&0
\end{eqnarray*}
whence the result follows.
\end{proof}

\begin{algorithm}\label{alg:L}
Given $A\in\Real^{n\times n}$ that is neutrally stable and
$C\in\Real^{m\times n}$, we obtain $L\in\Real^{n\times m}$ as
follows. Let $n_{1}\leq n$ be the number of eigenvalues of $A$ that
reside on the imaginary axis. Let $n_{2}:=n-n_{1}$.  If $n_{1}=0$,
then let $L:=0$; else construct $L$ according to the following steps.\\

\noindent
{\em Step 1:} Choose $U\in\Real^{n\times n_{1}}$ and $W\in\Real^{n\times n_{2}}$ satisfying
\begin{eqnarray*}
[U\ W]^{-1}A[U\ W]=
\left[
\begin{array}{cc}
F & 0\\
0 & G
\end{array}
\right]
\end{eqnarray*} 
where all the eigenvalues of $F\in\Real^{n_{1}\times n_{1}}$ have zero real parts. \\

\noindent
{\em Step 2:} Obtain $P\in\Real^{n_{1}\times n_{1}}$ from $F$ by \eqref{eqn:limit}. \\

\noindent
{\em Step 3:} Finally let $L:=UP^{-1}(CU)^{T}$.
\end{algorithm} 

\noindent
Below is our solution to Problem~I.

\begin{theorem}\label{thm:main}
Consider systems \eqref{eqn:system}. Let $u_{i}=Lz_{i}$ where
$L\in\Real^{n\times m}$ is constructed according to
Algorithm~\ref{alg:L} and $z_{i}$ is as in \eqref{eqn:z}. Then for all
network topologies described by connected $\Gamma$, solutions
$x_{i}(\cdot)$ for $i=1,\,2,\,\ldots,\,p$ synchronize to 
\begin{eqnarray*}
\bar{x}(t):=({r^{T}\otimes e^{At}})
\left[\!\!
\begin{array}{c}
x_{1}(0)\\
\vdots\\
x_{p}(0)
\end{array}
\!\!
\right]
\end{eqnarray*}
where $r\in\Real^{p}$ is such that $r^{T}\Gamma=0$ and $r^{T}\one=1$.
\end{theorem}

\begin{proof}
Let the variables that are not introduced here be defined as in
Algorithm~\ref{alg:L}.  Let $H:=CUP^{-1/2}$ and
$S:=P^{1/2}FP^{-1/2}$. Then $(H,\,S)$ is observable for $(C,\,A)$ is
detectable. Also, note that $S$ is skew-symmetric due to $PF+F^{T}P=0$.

We let $U^{\dagger}\in\Real^{n_{1}\times n}$ and $W^{\dagger}\in\Real^{n_{2}\times n}$ be such that
\begin{eqnarray*}
\left[
\begin{array}{c}
U^{\dagger}\\
W^{\dagger}
\end{array}
\right]=[U\ W]^{-1}\,.
\end{eqnarray*} 
Note then that $U^{\dagger}U=I_{n_{1}}$, $W^{\dagger}W=I_{n_{2}}$, $U^{\dagger}W=0$, and $W^{\dagger}U=0$.
Since $u_{i}=Lz_{i}$, we can combine \eqref{eqn:system} and \eqref{eqn:z} to obtain
\begin{eqnarray}\label{eqn:piece1}
\dot{x}_{i}=Ax_{i}+LC\sum_{j=1}^{p}\gamma_{ij}(x_{j}-x_{i})
\end{eqnarray}
Let now $\xi_{i}\in\Real^{n_{1}}$ and $\eta_{i}\in\Real^{n_{2}}$ be 
\begin{eqnarray}\label{eqn:piece2}
\left[
\begin{array}{c}
\xi_{i}\\
\eta_{i}\\
\end{array}
\right]:=
\left[
\begin{array}{cc}
P^{1/2}&0\\
0& I_{n_{2}}\\
\end{array}
\right]
\left[
\begin{array}{c}
U^{\dagger}\\
W^{\dagger}
\end{array}
\right]x_{i}
\end{eqnarray}
Combining \eqref{eqn:piece1} and \eqref{eqn:piece2} we can write
\begin{eqnarray}
\dot{\xi}_{i}&=&S\xi_{i}
+H^{T}H\sum_{j=1}^{p}\gamma_{ij}(\xi_{j}-\xi_{i})
+H^{T}CW\sum_{j=1}^{p}\gamma_{ij}(\eta_{j}-\eta_{i})\label{eqn:badem1}\\
\dot{\eta}_{i}&=&G\eta_{i}\,.\label{eqn:badem2}
\end{eqnarray}
Let $\Gamma$ be connected and $r\in\Real^{p}$ be such that
$r^{T}\Gamma=0$. Then define $\omega_{i}:\Real_{\geq
0}\to\Real^{n_{1}}$ as $\omega_{i}(t):=e^{-St}\xi_{i}(t)$ for
$i=1,\,2,\,\ldots,\,p$. Let $\dabilyu:=[\omega_{1}^{T}\ \omega_{2}^{T}\ \ldots\
\omega_{p}^{T}]^{T}$ and $\vi:=[\eta_{1}^{T}\ \eta_{2}^{T}\ \ldots\
\eta_{p}^{T}]^{T}$. Starting from \eqref{eqn:badem1} and
\eqref{eqn:badem2} we can write
\begin{eqnarray*}
\dot{\dabilyu}(t)&=&(\Gamma\otimes e^{-St}H^{T}He^{St})\dabilyu(t)
+(\Gamma\otimes e^{-St}H^{T}CWe^{Gt})\vi(0)\,.
\end{eqnarray*}
Thence
\begin{eqnarray}\label{eqn:integral}
\dabilyu(t)=\Phi(t,\,0)\dabilyu(0)
+\left[\int_{0}^{t}\Phi(t,\,\tau)(\Gamma\otimes e^{-S\tau}H^{T}CWe^{G\tau})d\tau\right]\vi(0)
\end{eqnarray}
where 
\begin{eqnarray*}
\Phi(t,\,\tau):={\rm exp}\left(
\int_{\tau}^{t}\left(\Gamma\otimes e^{-S\alpha}H^{T}He^{S\alpha}\right)d\alpha\right)
\end{eqnarray*}
is the state transition matrix \cite{antsaklis97}. From
Theorem~\ref{thm:lyap} we can deduce that $\Phi(t,\,\tau)$ is
uniformly bounded for all $t$ and $\tau$. Also, for any fixed $\tau$
we have $\lim_{t\to\infty}\Phi(t,\,\tau)=\one{r^{T}}\otimes
I_{n_{1}}$. Moreover, $e^{St}$ is uniformly bounded for all $t$, and
$e^{Gt}$ decays exponentially as $t\to\infty$ for $G$ is
Hurwitz. Therefore we can write
\begin{eqnarray*}
\lefteqn{\lim_{t\to\infty}\int_{0}^{t}\Phi(t,\,\tau)
(\Gamma\otimes e^{-S\tau}H^{T}CWe^{G\tau})d\tau}\\
&&\qquad\qquad\qquad=\int_{0}^{\infty}\left(\lim_{t\to\infty}\Phi(t,\,\tau)\right)
(\Gamma\otimes e^{-S\tau}H^{T}CWe^{G\tau})d\tau\\
&&\qquad\qquad\qquad=\int_{0}^{\infty}(\one{r^{T}}\otimes I_{n_{1}})(\Gamma\otimes e^{-S\tau}H^{T}CWe^{G\tau})d\tau\\
&&\qquad\qquad\qquad=0\,.
\end{eqnarray*}
Then, by \eqref{eqn:integral}, we can write
\begin{eqnarray*}
\lim_{t\to\infty}\dabilyu(t)=(\one{r^{T}}\otimes I_{n_{1}})\dabilyu(0)\,.
\end{eqnarray*}
Therefore solutions $\xi_{i}(\cdot)$ synchronize to $(r^{T}\otimes
e^{St})\dabilyu(0)$. Moreover, $\lim_{t\to\infty}\vi(t)=0$ for $G$ is
Hurwitz. Hence we can say that solutions $\eta_{i}(\cdot)$ synchronize to
$(r^{T}\otimes e^{Gt})\vi(0)$. As a result, solutions $x_{i}(\cdot)$ synchronize to
\begin{eqnarray*}
\left(r^{T}\otimes 
\left[UP^{-1/2}\ \ W\right]
\left[
\begin{array}{cc}
e^{St}&0\\0&e^{Gt}
\end{array}
\right]
\left[\begin{array}{c}P^{1/2}U^{\dagger}\\ W^{\dagger}\end{array}\right]
\right)
\left[
\begin{array}{c}
x_{1}(0)\\
\vdots\\
x_{p}(0)
\end{array}
\right]&&\\
=(r^{T}\otimes e^{At})
\left[
\begin{array}{c}
x_{1}(0)\\
\vdots\\
x_{p}(0)
\end{array}
\right]&&
\end{eqnarray*}
Hence the result.
\end{proof}

\section{Solution to Problem II}\label{sec:maindual}

This section, in which we provide a solution to Problem~II, follows
closely the previous one. We begin with the following algorithm.

\begin{algorithm}\label{alg:K}
Given $A\in\Real^{n\times n}$ that is neutrally stable and $B\in\Real^{n\times m}$,   
we obtain $K\in\Real^{m\times n}$ as follows. Let $n_{1}\leq n$ be the number of 
eigenvalues of $A$ that reside on the imaginary axis. Let $n_{2}:=n-n_{1}$. 
If $n_{1}=0$, then let $K:=0$; else
construct $K$ according to the following steps. \\

\noindent
{\em Step 1:} Choose $U\in\Real^{n\times n_{1}}$ and $W\in\Real^{n\times n_{2}}$ satisfying
\begin{eqnarray*}
[U\ W]^{-1}A[U\ W]=
\left[
\begin{array}{cc}
F & 0\\
0 & G
\end{array}
\right]
\end{eqnarray*} 
where all the eigenvalues of $F\in\Real^{n_{1}\times n_{1}}$ have zero
real parts. Let $U^{\dagger}\in\Real^{n_{1}\times n}$ and
$W^{\dagger}\in\Real^{n_{2}\times n}$ be such that
\begin{eqnarray*}
\left[
\begin{array}{c}
U^{\dagger}\\
W^{\dagger}
\end{array}
\right]=[U\ W]^{-1}\,.
\end{eqnarray*}\\

\noindent
{\em Step 2:} Obtain $P\in\Real^{n_{1}\times n_{1}}$ from $F$ by \eqref{eqn:limit}. \\

\noindent
{\em Step 3:} Finally let $K:=(U^{\dagger}B)^{T}PU^{\dagger}$.
\end{algorithm} 

\noindent
Below is our solution to Problem~II.

\begin{theorem}\label{thm:maindual}
Consider systems \eqref{eqn:dualsystem}. Let $u_{i}=Kz_{i}$ where
$K\in\Real^{m\times n}$ is constructed according to
Algorithm~\ref{alg:K} and $z_{i}$ is as in \eqref{eqn:zdual}. Then for
all network topologies described by connected $\Gamma$, solutions
$x_{i}(\cdot)$ for $i=1,\,2,\,\ldots,\,p$ synchronize to 
\begin{eqnarray*}
\bar{x}(t):=({r^{T}\otimes e^{At}})
\left[\!\!
\begin{array}{c}
x_{1}(0)\\
\vdots\\
x_{p}(0)
\end{array}
\!\!
\right]
\end{eqnarray*}
where $r\in\Real^{p}$ is such that $r^{T}\Gamma=0$ and $r^{T}\one=1$.
\end{theorem}

\begin{proof}
Let the variables that are not introduced here be defined as in
Algorithm~\ref{alg:K}.  Let $H:=(P^{1/2}U^{\dagger}B)^{T}$ and
$S:=P^{1/2}FP^{-1/2}$. Then $(S,\,H^{T})$ is controllable for $(A,\,B)$ is
stabilizable. Also, note that $S$ is skew-symmetric due to $PF+F^{T}P=0$.

Since $u_{i}=Kz_{i}$, we can combine \eqref{eqn:dualsystem} and \eqref{eqn:zdual} to obtain
\begin{eqnarray}\label{eqn:piece1dual}
\dot{x}_{i}=Ax_{i}+BK\sum_{j=1}^{p}\gamma_{ij}(x_{j}-x_{i})
\end{eqnarray}
Let now $\xi_{i}\in\Real^{n_{1}}$ and $\eta_{i}\in\Real^{n_{2}}$ be 
\begin{eqnarray}\label{eqn:piece2dual}
\left[
\begin{array}{c}
\xi_{i}\\
\eta_{i}\\
\end{array}
\right]:=
\left[
\begin{array}{cc}
P^{1/2}&0\\
0& I_{n_{2}}\\
\end{array}
\right]
\left[
\begin{array}{c}
U^{\dagger}\\
W^{\dagger}
\end{array}
\right]x_{i}
\end{eqnarray}
Combining \eqref{eqn:piece1dual} and \eqref{eqn:piece2dual} we can write
\begin{eqnarray}
\dot{\xi}_{i}&=&S\xi_{i}
+H^{T}H\sum_{j=1}^{p}\gamma_{ij}(\xi_{j}-\xi_{i})\label{eqn:badem1dual}\\
\dot{\eta}_{i}&=&G\eta_{i}+W^{\dagger}BH\sum_{j=1}^{p}\gamma_{ij}(\xi_{j}-\xi_{i})\,.
\label{eqn:badem2dual}
\end{eqnarray}
Looking at \eqref{eqn:badem1dual}, by Theorem~\ref{thm:lyap}, we 
assert that solutions $\xi_{i}(\cdot)$ synchronize to
\begin{eqnarray*}
({r^{T}\otimes e^{St}})
\left[\!\!
\begin{array}{c}
\xi_{1}(0)\\
\vdots\\
\xi_{p}(0)
\end{array}
\!\!
\right]
\end{eqnarray*}
Now observe that $|\xi_{j}(t)-\xi_{i}(t)|\to 0$ exponentially as $t\to
\infty$ for all $(i,\,j)$ pairs.  Also recall that $G$ is
Hurwitz. From \eqref{eqn:badem2dual} we can therefore deduce by
input-to-state stability (ISS) arguments \cite{khalil96} that
$\eta_{i}(t)\to 0$ as $t\to\infty$ for $i=1,\,2,\,\ldots,\,p$. The remainder of the 
proof is same as that of proof of Theorem~\ref{thm:main}. 
\end{proof}

\section{Conclusion}
Let us now briefly discuss the generality of the assumptions in the
paper.  For linear time-invariant case with identical individual
system dynamics, it should be evident that detectability
(stabilizability) assumption is indispensable for
synchronization. Regarding the neutral stability condition, it would
be of great interest to study the synchronization of unstable
systems. However, when neutral stability assumption on individual
systems is relinquished, mere connectedness of the network should
generally not be sufficient for individual systems to synchronize. The
reason is that, due to unstable dynamics, the trajectories will tend
to drift apart from each other when there is no (or very little)
coupling. The coupling strength therefore should be above some
threshold to overcome that tendency, which requires a stronger (than
connectedness) condition on the network topology.  

\bibliographystyle{plain}         
\bibliography{references}            

\end{document}

%% file: macros.tex
\usepackage[final]{graphicx}   
\usepackage{psfrag}   
\usepackage{amsmath,amsfonts,amssymb,amsxtra}

\newcommand {\Real}{\ensuremath{{\mathbb{R}}}}
\newcommand {\Natural}{\ensuremath{{\mathbb{N}}}}

\newcommand{\A}{\ensuremath{\mathcal A}}

\newcommand{\X}{\ensuremath{\mathcal X}}

\newcommand{\N}{\ensuremath{\mathcal N}}
\newcommand{\W}{\ensuremath{\mathcal W}}

\newcommand{\vi}{\ensuremath{{\mathbf{v}}}}
\newcommand{\ex}{\ensuremath{{\mathbf{x}}}}

\newcommand{\one}{\ensuremath{{\mathbf{1}}}}

\newcommand{\dabilyu}{\ensuremath{{\mathbf{w}}}}

\newtheorem{theorem}{Theorem}
\newtheorem{corollary}{Corollary}

\newtheorem{fact}{Fact}
\newtheorem{algorithm}{Algorithm}

\newenvironment{proof}{\noindent {\bf Proof.}}{\hfill \hspace*{1pt}\hfill$\blacksquare$}